\documentclass[12pt]{amsart}
\usepackage{amsmath}
\usepackage{epsfig,amssymb,epic,eepic,color}
\usepackage[all]{xy}
 \usepackage[T1]{fontenc}  
 \usepackage{hyperref}    
 \usepackage{textcomp}  
 
\usepackage{xcolor}

\newcommand{\be}{\begin{enumerate}}
\newcommand{\ee}{\end{enumerate}}
\newcommand{\bi}{\begin{itemize}}
\newcommand{\ei}{\end{itemize}}

\def\R{\mathbb{R}}

\def\Z{\mathbb{Z}}

\def\om{\omega}

\def\ga{\gamma}    
\def\Ga{\Gamma}

\def\al{\alpha}
\def\be{\beta}
 
\def\de{\delta}
\def\De{\Delta}

\def\la{\lambda}
\def\La{\Lambda}

\def\Si{\Sigma}

\def\ep{\varepsilon}

\def\nd{\noindent}
\def\bull{\hfill$\Box$\\}

\begin{document}
\vskip 1cm

\begin{center}
{\sc  A geometric Morse-Novikov complex with infinite series coefficients
\medskip

  Fran\c cois Laudenbach \& Carlos Moraga Ferr\'andiz

}

\vskip .5cm

\end{center}
\title{}
\address{Laboratoire de Math\'ematiques Jean Leray,  UMR 6629 du CNRS, Facult\'e des Sciences et Techniques,Universit\'e de Nantes, 2, rue de la Houssini\`ere, F-44322 Nantes cedex 3, France.}
\email{francois.laudenbach@univ-nantes.fr}
 
 \address{36, av. Camille Gu\'erin, 44000 Nantes, France}
\email{crlsmrgf@gmail.com}



\vskip .5cm

\begin{itemize}
{\small 
\item[] {\sc Abstract.} Let $M$ be a closed $n$-dimensional manifold, $n>2$, whose 
 first real cohomology group $H^1(M;\R)$ is non-zero.
We present a general method for 
constructing a Morse 1-form $\al$ on $M$, closed but non-exact,
 and a pseudo-gradient $X$ such that the differential 
$\partial^X$ of the Novikov complex of the pair $(\al, X)$ has at least one incidence coefficient which is an infinite series. 
This is an application of our previous study of the homoclinic bifurcation of pseudo-gradients of multivalued Morse 
functions.
}\\
\end{itemize}


\thispagestyle{empty}
\vskip .5cm



\section{Introduction}
\medskip
Let us consider an $n$-dimensional closed smooth manifold $M$ with a non-zero de Rham cohomology in degree one.
Let $u\in H^1_{dR}(M)$, $u\neq 0$,  and let $\al$ be a differential closed form in the class $u$. By Thom's
transversality theorem with constraints \cite{thom-jet,thom56}, generically in the class $u$, the local primitives
of a closed 1-form $\al$ are Morse functions (such an $\al$ is named a Morse 1-form).
 The finite set $Z(\al)$ made of its zeroes is graded by the Morse index: $Z(\al)= \cup_{k=0}^n  Z_k(\al)$.
 
 Taking the point of view that a Morse 1-form is nothing but a multivalued function (that is, a function well-defined up to an 
 additive constant), S.P. Novikov developed a theory analogous to  Morse theory \cite{nov1,nov2} (today 
 referred to as {\it Morse-Novikov theory}) which leads to some
 homology denoted by $H_*(M; u)$. Here, we summarize it with the point of view adopted by J.-C. Sikorav \cite{sikorav}.
 
 The class $u$ is seen as the morphism $\pi_1(M)\to \R$ defined by 
 $u(\ga)=\int_\ga\al$ for every oriented loop $\ga$.
 Let $\tilde u: \Z\pi_1(M)\to\overline{\R}$ be the associate valuation: for $\la\in\Z\pi_1(M)$, that is a finite sum
 $\sum_i n_i g_i$ with $n_i\in \Z^*$ and $g_i\in\pi_1(M)$, one defines $\tilde u(\la)= \max u(g_i)$.
 The {\it universal} Novikov ring $\La_u$ is the completion of 
 $ \Z\pi_1(M)$ for $\tilde u(\la)\to -\infty$. It is worth noticing that we get the same completion by replacing 
 $\pi_1(M)$ with the {\it fundamental groupoid} of $M$. This choice is made in what follows.
 
 The Novikov complex $N_*(\al,\partial^X)$ is the free $\La_u$-module based on 
 $Z_*(\al)$; the differential $\partial ^X$ is detailed right below.
 For defining $\partial ^X$, one chooses some (descending) pseudo-gradient $X$ {\it adapted} to
 $\al$, meaning that $X$ fulfils the following:
 
 1) $\langle \al(x), X(x)\rangle<0$ for every $x\in M$ apart from $Z(\al)$;
 
 2) for every $p\in Z_k(\al)$, the field $-X$ is the Euclidean gradient of a local primitive $h_p$
 of $\al$ in Morse coordinates $(x_1, \ldots, x_n)$ about $p$ 
 where  $h_p$ reads $$h_p(x_1,\ldots, x_n)= -x_1^2-\ldots-x_k^2+x_{k+1}^2\ldots +x_n^2.$$

 The zeroes of $X$, which coincide with the zeroes of $\al$, are {\it hyperbolic}. Therefore, for each $p\in Z_k(\al)$, 
 we have a stable manifold $W^s(p,X)$ of dimension $n-k$ and 
 an unstable manifold $W^u(p,X)$ of dimension $k$.
  According to  Kupka-Smale's theorem (see \cite{palis}), one can approximate $X$
 so that all stable and unstable manifolds are mutually transverse. For $p\in Z_k(\al)$, the value of $\partial ^X$ in $p$
 reads:
 $$\partial^Xp=\sum_{q\in Z_{k-1}(\al)}\langle p,q\rangle^X q
$$
where the {\it incidence coefficient} $\langle p,q\rangle^X$ belongs to the Novikov ring $\La_u$.

As far as we know, there are no published    
example\footnote{\label{note} A 3-dimensional example where $u$ is a rational class has been given by  A. Pajitnov in
  \cite[Section 3]{pajitnov}.} 
where some of the incidence coefficients are infinite series as the 
Novikov ring allows it. By contrast, it is known that if $f:M\to \R$ is a Morse function, then for every $c\gg 0$, the 
differential of the Novikov complex $N_*(\al+c\,df,X)$ has all its incidence coefficients in $\Z\pi_1(M)$. Here, $X$ denotes
a pseudo-gradient adapted to $\al+c\,df$ which is $C^1$-close to a pseudo-gradient adapted to $f$ \cite[Lemma 3.7]{harvey}. 
When $n>2$, we present a general method for constructing some pairs $(\al, X)$ such that the differential $\partial ^X$
of $N_*(\al)$ have some incidence coefficients which are infinite series. Our method is based on our study 
 \cite{l-m} of the 
homoclinic bifurcations of the pseudo-gradient $X$.\\

\nd {\sc Theorem.} {\it Let $M$ be an $n$-dimensional closed manifold with $n>2$. Let $u\in H^1_{dR}(M)$ be a non-zero
 cohomology class and let $\al$ be a Morse 1-form in the class $u$ without centers, that is, with 
no zeroes of extremal index. 
 For  $1\leq k\leq n-2$, let $\al'$ be any Morse 1-form obtained from $\al$ by creating the birth of a pair of zeroes $(q,p)$
 of respective indices $(k, k+1)$. Then, it is possible to build a pseudo-gradient  $X'$ adapted to $\al'$ such that the incidence coefficient $\langle p,q\rangle^{X'}$ is an infinite series in the Novikov ring $\La_u$.}\\
 
 \nd{\sc Remark.} Observe that, when 
 $\al$ is a non-exact Morse 1-form,
  it is possible to cancel all zeroes of 
 extremal indices \cite[Remarque section 1.2]{arnoux}. 
 Therefore, the assumption of having no centers is irrelevant.
 
 \section{Proof}
 \subsection{\sc Making a tube of periodic orbits.}
 Let $X$ be a Kupka-Smale pseudo-gradient adapted to $\al$. 
 After the no-center assumption, almost no orbit of $X$ goes to (or comes from) one of its zeroes. A $C^0$-approximation
 of $X$ near  an accumulation point of such an orbit allows us to create not only one periodic orbit but a tube
 $T\cong D^{n-1}\times S^1$ made of periodic orbits $\{pt\}\times S^1$ (still name $X$ the pseudo-gradient after this approximation). 
  The (free) homotopy class of these periodic orbits is denoted  by $g$. 
  Note $u(g)<0$.\\

 \subsection{\sc Birth of a pair of zeroes in cancelling position.} Let $B\subset T$, $B \cong D^{n-1}\times D^1$, 
  be a polydisc  bi-foliated 
  by leaves of $\al$ and orbits of $X$. Let $h$ be a primitive of $\al\vert_{ B}$. 
  Let us insert  a {\it birth model} into $B$ to create a pair $(p,q)$ of critical points of respective indices $k+1$ and 
  $k$ (see \cite[chap. III]{cerf} and also \cite{whitney}). We impose this one-parameter 
  deformation starting from $h$ to be supported in $B$.
  Denote $h'$ the function ending this deformation; set $\al'= \al +dh'-dh$. 
   Let $X_1$ be a pseudo-gradient adapted to $\al'$ near $p$ and $q$ and coinciding with $X$ outside $B$.
    There is exactly
   one orbit of $X_1$ descending from $p$ to $q$ inside $B$. 
  However, 
  the pseudo-gradient $X_1$ is not convenient for our theorem as some extra
  connecting orbit exists in $T$
  from $p$ to $q$. Thus, we change $X_1$ into $X'$ as follows:
  make $\partial T$ an attractor of $X'\vert_{ T}$ (while keeping 
   $X'\vert_{\partial T}=X_{1}\vert_{\partial T}$)
  and one of the periodic orbits in $T\smallsetminus B$ a repeller
  in order that any orbit of $X'$ emanating from the bottom of $B$, namely 
  $\partial^-B:=\{h'=h'(q)-\ep\}$, $0<\ep\ll 1$,
   never returns to $B$ in positive time. The outcome of this construction is that the Novikov complex
   $N_*(\al', X')$ fulfils the following: $\partial ^{X'}q=0 $ and $\partial ^{X'}p=q$ for suitable orientations 
   of the unstable manifolds. Moreover, 
   as the boundary of $T$ is invariant by $X'$, no 
   zero of $\al'$ other than $p$ has a connecting orbit to $q$; and no zero is connected to $p$.
   
   Let $\mathcal S_g$ be the stratum of pseudo-gradients adapted to $\al'$ which are equal to
     $X'$  outside $T$ and have one 
    homoclinic connecting orbit  from $p$ to $p$ in $T$ in the class $g$ with the minimal defect of transversality.
    This is a codimension-one stratum which is canonically co-oriented.
     \\
    \subsection{\sc Notation.} Denote by $\partial^+B:= \{h'=h'(p)+\ep\}$ the top of $B$. It contains the so-called 
   {\it belt sphere} 
   $\Si^+$ of $p$, that is, $\Si ^+ : = W^s(p, X')\cap \partial^+B$. This sphere is $(n-k-2)$-dimensional and co-oriented by 
   the chosen orientation of $W^u(p,X')$. Moreover, it is the 
   boundary of the closed $(n-k-1)$-disc
   $\De^+$ whose interior is made of points in $\partial^+B$ whose positive $X'$-flow goes to $q$. Similarly, set
   $\partial^-B: =\{h'=h'(q)-\ep\}$. Denote by $\De^-$ the closure of  $W^u(p,X')\cap \partial^-B$ (a $k$-disc) and 
   set $\Si^-:= W^u(q,X')\cap \partial^-B$. 
   The sphere $\Si^-$ bounds $\De^-$. 
     
    Let $\mathcal L^\pm$ denote the leaf of $\al'$ in $T$
    which contains $\partial^\pm B$. 
     The ball $\partial ^+B$ shall be regarded as the union of a normal   tube 
     $N(\Si^+)$ to $\Si^+$ and a normal 
     tube   
     $N(\De^+)$ to $\De^+$ in $\mathcal L^+$. Moreover, there exists a larger tube $N'(\Si^+)$ which 
     contains all fibres of $N(\De^+)$ which meet $N(\Si^+)$.

   Let $\rho: \mathcal L^-\to\mathcal L^+$ be the {\it first} holonomy diffeomorphism 
   of the positive flow of $X'$.  
   Denote by $Desc: \mathcal L^+\smallsetminus \De^+\to \mathcal L^-\smallsetminus \De^-$ 
   the {\it gradient descent} by 
   $X'$. 
   Set $D:= \rho(\De^-)$; 
   it is a $k$-disc in $\mathcal L^+$ disjoint from $\partial^+B$ whose boundary lies in $W^u(q,X')$.
   
   Finally, 
   let $\mu_a$ be the 
   fibre in $a\in \Si^+$ of the tube $N'(\Si^+)$.
   By descent, the $(k+1)$-disc $\mu_a$ descends to $\Ga^-_a$, a cylinder over $\De^- $, pinched along $\partial \De^-$. 
   The cylinder $\Ga^-_a$ is foliated by {\it rays} which are traced by the normal fibres to $\De^-$ in $\mathcal L^-$.
    The holonomy diffeomorphism
    $\rho$ carries $\Ga^-_a$ and its foliation to a foliated pinched cylinder $\Ga^+_a$ whose leaves are still named rays.

   \subsection{\sc Some isolated self-slide of $p$; homoclinic bifurcation.} 
   In this setup, a {\it positive self-slide} of $p$ consists in making an isotopy of 
    $D$ in $\mathcal L^+$, relative to $\partial D$,
   which crosses $\Si^+$ exactly once and positively with respect to its co-orientation \cite[Remark 3.2]{l-m}. 
   This isotopy lifts to a deformation $(X'_s)_{s\in [0,1]}$ of $X'_0=X'$.
   At the time of crossing (let us say $s= \frac 12)$, the homoclinic 
   connecting orbit from $p$ to $p$ lies in 
    the homotopy class $g$. The path
   $(X'_s)_{s}$ is said to cross the codimension-one stratum $\mathcal S_g$ positively 
   in the space of pseudo-gradients adapted to $\al'$.
   
   Up to homotopy, 
   the path $(X'_s)_s$ is well-defined by the following data: 
   \begin{itemize}
   \item[--] a simple smooth path  $\ga:[0, \frac 12]\to \mathcal L^+$, starting non-tangentially
   to $D$ 
   from any point 
   $\ga(0)\in int(D)$, 
   ending 
   in $a^+:=\ga(\frac 12)$ normally to $\Si^+$ 
    and which, for every $s\in (0,\frac 12)$, avoids  
   both of $W^u(p,X')\cup \De^+$ and the repelling periodic orbit.
    \item[--] a transverse framing $\tau=(\tau_1, \ldots,\tau_{n-2})$ to $\ga$ in $\mathcal L^+$ such that
    $\tau_1^k:=\tau_1\wedge\cdots\wedge \tau_k$
   is tangent to $D$ in $\ga(0)$ and such that $(\dot\ga(\frac 12), \tau_1,\ldots, \tau_{n-2})$ 
   is an orthogonal framing 
   positively normal to $\Si^+$ in $a^+$.
   \end{itemize}
    Note that any point  $a^+\in\Si^+$ is reachable from $D$ by such a path. In what follows, we will specify 
     $\ga$ more precisely near $s\in\{0, \frac 12\}$. For further use, we introduce $a^-: = \rho^{-1}(\ga(0))$, 
     a point of $\De^-$.
      Let $\tilde\Si^-:= W^u(p,X')\cap\{h'=h'(p)-\ep\}$ be the attaching $k$-sphere associated
   with $p$ and denote by  $\tilde a ^-\in \tilde\Si^-$ the point which descends to $a^-$ through the flow of $X'\vert_B$. 
   
     Let $Y$ be a vector field on $\mathcal L^+$, tangent to $\De^+$,
   whose support $S$ is contained in a 
   neighbourhood of $\ga([0, \frac 12])\cup \partial^+B$ and
    whose flow is noted $Y^s$, $s\in \R$. 
    The triple $(\ga,\tau,Y)$ may be chosen so that the following holds.
   \begin{itemize}
  \item[(i)] For every $s\in [0,\frac 12]$, we have $Y^s(\ga(0))= \ga(s)$ and $\tau_1^k$ is tangent
   to $Y^s(D)$ in $\ga(s)$; moreover, we impose the vector $\dot\ga(\frac12)$ to be tangent to $\De^+$ pointing 
   inwards.
   \item[(ii)]  Denoting by $\mathcal T$  the closure of $\mathop{\cup}_{s=0}^{+\infty}Y^s(D)$, named the \emph{tongue}, 
   the intersection $\mathcal T\cap \partial^+B$ is a corned 
   disc $\hat \nu_{a^+}$ contained in 
   the fibre $\mu_{a^+}$ of the tube
   $N'(\Si^+)$. Denote by $\nu_{a^+}$ the corresponding fibre of $N(\Si^+)$.
   
  \item[(iii)] For every $s\in[\frac 12, +\infty)$, the  intersection $Y^s(D)\cap \partial ^+B$ is one fibre of $N(\De^+)$
  which meets
  $N(\Si^+)$. The  distance of this fibre to $\Si^+$ is noted $r(s)$; set $\eta= \max r(s)$. A suitable choice of $Y$  
  makes it as small as desired.
  
  \item[(iv)] The vector field $Y$ is required to be tangent to the foliation of $\Ga^+_{a^+}$ by rays. This requirement forces 
  the path $\ga$ to be tangent to a ray of $\Ga^+_{a^+}$ and to follow it entirely.
  
 \item[(v)] The support $S$ satisfies the following condition: $(S\setminus \partial^+ B)$ is contained in a small
 neighbourhood of $\ga([ 0, \frac 12])$.
  \end{itemize}
  The last requirement needs some preparation.  By (i), the vector $\dot \ga(\frac 12)$ is tangent to $\De^+$ and points
     in the direction  marked by the point $\De^+\cap \partial\nu_{a^+}$. 
     This vector determines a co-oriented \emph{equator}
    $E$ in the sphere $\partial \nu_{a^+}$ (see \cite[subsection 2.4]{l-m}).
    Recall that the Morse model of $p$ produces 
    a canonical diffeomorphism  from the meridian sphere $\partial\nu_{a^+}$ to $\tilde\Si^-$. It  maps $E$ to
     an equator $\tilde E\subset \tilde\Si^-$ which avoids the connecting orbit from $p$ to $q$ 
      by the choice of $\dot \ga(\frac 12)$. Thus, 
    $\tilde E$  descends through the flow of $X'\vert_B$
    to a $(k-1)$-sphere $E^-$ contained in the interior of $ \De^-$. This $E^-$ makes a decomposition 
    $\De^-=\de^-\mathop{\cup}_{E^-}\de^+$ where $\de^-$ is a $k$-disc and $\de^+$ 
    is diffeomorphic to $ S^{k-1}\times[0,1]$. The signs of $\de^\pm$  
    reflect the co-orientation of $E$ in $\partial \nu_{a^+}$, 
    and hence, of $E^-$ in $\De^-$.
     By moving $\ga(0)$ in $D$, we can fulfill the last requirement below. 
    \begin{itemize}
    \item[(vi)] The point $a^-$ lies in $\de^+$. Moreover, we make 
     the support $S$ of $Y$ disjoint from $\rho(\de^-)$.
     \end{itemize}
     \nd  Finally, the isotopy $(Y^s)_{s\in [0,1]}$ 
  lifts to the desired one-parameter deformation  $(X'_s)$ of $X'$  realizing a self-slide of $p$ in the homotopy class $g$.\\
  
  \nd{\sc Lemma 1.} {\it With the above requirements {\rm (i-vi)}, 
  the self-slide at $s=\frac 12$ is \emph{isolated}. More precisely, for
   $s\neq \frac 12$, the pseudo-gradient $X'_s$ has no homoclinic orbit.}\\
   
   \nd{\sc Proof.} The dynamics of $\rho$ which moves $\partial D$ away from $S$ makes clear that 
   $W^u(q, X'_s)$ never meet $\De^+$. This implies that no $s\in [0,1]$ is a time of self-slide of $q$.  
   
   Concerning the self-slides of $p$,  we argue as follows. By (ii),    the  corned disc $\hat\nu_{a^+}$ 
    descends to 
    a  corned cylinder $C^-\subset \Ga^-_{a^+}$, also based on $\De^-$ and 
    pinched along $\partial \De^-$.
    Then, the holonomy diffeomorphism
    $\rho$ carries $C^-$ to a pinched corned cylinder 
    $C^+ \subset \Ga^+_{a^+}$.   
    By (iv-v), the intersection $C^+ \cap S$ is 
    foliated by (truncated) rays of $C^+$  starting from $D\cap S$. 
        In particular, $S\cap \mathcal T$ may be viewed as a \emph{sub-tongue} of $\mathcal T$.

     We are now ready to prove Lemma 1. We limit ourselves to prove that $X'_s$ is never a self-slide of $p$ in the 
     class $g^2$ when $s\neq \frac 12$, that is:  
     $$(*)\quad \forall s\neq \frac 12,\quad [(Y^s\rho)\circ Desc\circ_1(Y^s\rho)(\De^-)]\cap \Si^+=\emptyset \ .$$
     Here, the operator $\circ_1$ stands for $Desc$ where it is defined, that is 
     in $\mathcal L^+\smallsetminus \De^+$.
     The proof concerning the classes $g^k$, $k>2$, is similar with the same choices (i) - (vi). 
     
     When $s<\frac 12$ 
     and as long as  the disc $Y^s(D)$ does not touch the ball $\partial^+ B$, the dynamics of $\rho\circ Desc$
     implies $(*)$. When $Y^s(D)\cap N(\Si^+)$ is non-empty, 
     we define $D_s^-:= Desc\circ_1Y^s(D)$. 
     When $s<\frac 12$,  one checks that $D_s^-$ is contained in the union of 
     thin tube $Desc(S\smallsetminus \partial^+B)$ and the part 
      of the cylinder $C^-$ {\it over} $\de^-$ (in the sense of the rays). 
     Thus by (vi),
     $\rho(D^-_s)$ does not approach the support  $S$ of $Y$ 
     and $(*)$ holds for every $s\in [0, \frac 12)$.
     
     For $s>\frac 12$, 
     $D^-_s$ is now an annulus whose boundary is made of two $(k-1)$-spheres, one of which is
     $Desc(\partial D)$ and the other is $\partial \De^-$. 
      By \cite[Lemma 3.8]{l-m}, 
     when $s$ goes to $\frac 12$ from above, it approaches any 
     compact domain given in the interior of the annulus $\de^+$ in the $C^1$-topology. As a consequence of (iv), 
     when $s-\frac 12$ is small enough the intersection 
     of $\rho(D^-_s)$ with $S$ is the graph $G_{u_s}$ 
     in the cylinder 
     $C^+$ of some positive 
     function $u_s: \rho(\de^+)\cap S\to \R_{>0}$.
      Since the flow of $Y$ preserves the order on each orbit, $[(Y^s\rho)(D^-_s)]\cap S$ 
     (which lies in the tongue) is
     in front position with respect to $[Y^{1/2}(D)]\cap S$ on each flow line of $Y$.  This claim holds
      true when $s>\frac 12$ is very close to $\frac 12$ and it remains true up to 1 by the flow rule.
     Therefore, $[(Y^s\rho)(D^-_s)]\cap S$  is still beyond $[Y^{1/2}(D)]\cap S$ in the tongue, and hence, it
      cannot meet $\Si^+$.  
      ${}$ \hfill \bull

   \subsection{\sc The character function.} This function
    is a natural continuous function $\chi:\mathcal S_g\to\R$ which takes 
   positive \emph{and} negative values on each connected component of $\mathcal S_g$ \cite[Theorem 1.1]{l-m}. 
   For a pseudo-gradient $Y\in \mathcal S_g$, the value $\chi(Y)$ 
    only depends on the 1-jet of the holonomy $\frak h$
   of $Y$ along the unique homoclinic orbit $\ell_Y$ from a zero of $Y$ to 
   itself in the homotopy class $g$.
   The  construction and properties of $\chi(Y)$ are given
   in \cite[Section 2]{l-m}. 
   Denote by $\mathcal S_g^+$ (resp. $\mathcal S_g^-$)
   the non-empty open set in $\mathcal S_g$ where $\chi>0$ (resp. $\chi<0$). When crossing $\mathcal S_g$ transversely,
   the change of the differential in 
   the Morse-Novikov complex is different according to  the crossing: it can be positive or negative;
   through $\mathcal S_g^+$ or $\mathcal S_g^-$. A positive crossing through  $\mathcal S_g^+$ implies a change
   by multpiplication with the infinite series $1+g+g^2+ \cdots$.

   In our setting, the considered pseudo-gradient is $X'_{1/2}$ and the involved homoclinic orbit $\ell$
   is the  orbit  of $X'_{1/2}$ which leaves the box  $B$ at $a^-$ and enters $B$ at $a^+$.
   We wish $\chi(X'_{1/2})>0$. By the construction, the character reads 
   $\chi(X'_{1/2})= \la \om_+ +\om_-$ with $\la>0$ and $\om_\pm\in [-1,+1]$. Each term\footnote{The terms 
   $\om^+$ and $\om^-$ are named the {\it latitudes} in \cite{l-m} and the factor $\la$ is called the {\it holonomy factor}.} 
   in this formula
   depends on $\frak h$ only.
   The term $\om_-$ is positive by the choice of $\dot\ga(\frac 12)$
   (which fixes the decomposition of $\De^-= \de^-\cup_{E^-}\de^+$) and  by taking $a^-\in \de^+$.
      The term $\om_+$ is not estimated but the factor $\la$ may be made as small as desired
   by the choice of the holonomy of $Y$
    along $\ga([0,\frac 12])$ once $\ga(0)= \rho(a^-)$ and $\ga(\frac 12)=a^+$ are fixed.  And this choice
    is free. 
   As a consequence,  
   by \cite[Theorem 3.6]{l-m}, we have:
   $$
  (*)\quad  \left\{
  \begin{array}{l}
  \partial^{X'_1} q=0,\\
  \partial^{X'_1} p= (1+g+g^2+g^3+\cdots)\, (\partial^{X'_0} p)=  (1+g+g^2+g^3+\cdots) \,q\,.
  \end{array}\right.
   $$
   
   Some  comment about $(*)$ is necessary here.
   In the theorem which is referred to, $(*)$ is stated only for the germ of $(X'_s)_s$ at $s=\frac 12$ 
   with truncation of the incidence coefficients. The support of a representative of the germ depends  
   on the order of the truncation, the higher the order of truncation 
   the smaller the domain of the representative. 
   Here,  we are in a very particular situation: by Lemma 1
    the self-slide is isolated. 
    As a consequence, the truncation in question is not needed.
    This finishes the proof of our theorem.

\end{document}